\documentclass[11pt, a4paper]{amsart}
\usepackage{a4wide}
\usepackage{mathrsfs, amssymb,amsmath,url,amsthm}

\title{Permutations on the random permutation}
\date{\today}

\author[J.~Linman]{Julie Linman}
	\address{Department of Mathematics\\ University of Colorado\\ Boulder, CO 80309-0395, USA}
	\email{julie.linman@colorado.edu}
\author[M.~Pinsker]{Michael Pinsker}
    \address{\'{E}quipe de Logique Math\'{e}matique\\ Universit\'{e} Diderot -- Paris 7\\
	UFR de Math\'{e}matiques\\
	75205 Paris Cedex 13, France}
    \email{marula@gmx.at}
    \urladdr{http://dmg.tuwien.ac.at/pinsker/}
\thanks{The research of Michael Pinsker has been funded through project I836-N23 of the  Austrian Science Fund (FWF)}

\theoremstyle{plain}

    \newtheorem{thm}{Theorem}[section]

    \newtheorem{lem}[thm]{Lemma}

    \newtheorem{prop}[thm]{Proposition}

    \newtheorem{cor}[thm]{Corollary}

    \newtheorem{fact}[thm]{Fact}
    
    \newtheorem{prob}[thm]{Problem}

\theoremstyle{definition}

    \newtheorem{defn}[thm]{Definition}

\theoremstyle{remark}

\usepackage{tikz}
\usetikzlibrary{arrows,shapes,positioning,trees,mindmap,backgrounds}

\tikzstyle{subgroup}=[rounded rectangle,draw=black,scale=.75]
\tikzstyle{JI}=[rounded rectangle,draw=black,thick,scale=.75]
\tikzstyle{line}=[draw]

\DeclareMathOperator{\id}{id}
\DeclareMathOperator{\st}{St}
\DeclareMathOperator{\tw}{Tw}
\DeclareMathOperator{\sw}{sw}

\DeclareMathOperator{\btw}{Btw}
\DeclareMathOperator{\cyc}{Cyc}
\DeclareMathOperator{\sepa}{Sep}
 
\DeclareMathOperator{\Aut}{Aut}
\DeclareMathOperator{\Sym}{Sym}

\DeclareMathOperator{\JI}{JI}

\newcommand{\ignore}[1]{}

\newcommand{\cl}[1]{\langle #1 \rangle}

\newcommand{\C}{{\mathscr C}}

\newcommand{\F}{{\mathscr F}}

\newcommand{\Q}{{\mathbb Q}}

\newcommand{\To}{\rightarrow}

\newcommand{\G}{{\mathscr G}}
\renewcommand{\H}{{\mathscr H}}

\newcommand{\rest}{{\upharpoonright}}

\renewcommand{\L}{{\mathfrak L}}

\newcommand{\oo}{<_1}
\newcommand{\ot}{<_2}
\DeclareMathOperator{\iiii}{int}
\DeclareMathOperator{\rev}{rev}
\DeclareMathOperator{\up}{Up}
\DeclareMathOperator{\dow}{Do}
\renewcommand{\int}{\iiii}

\newcommand{\la}{\langle}
\newcommand{\ra}{\rangle}
\newcommand{\mix}[2]{\binom{#1}{#2}}

\begin{document}

\begin{abstract}
The \emph{random permutation} is the Fra\"{i}ss\'{e} limit of the class of finite structures with two linear orders. 
Answering a problem stated by Peter Cameron in 2002, we use a recent Ramsey-theoretic technique to show that there exist precisely 39 closed supergroups of the automorphism group of the random permutation, and thereby expose all symmetries of this structure. Equivalently, we classify all structures which have a first-order definition in the random permutation.

\end{abstract}
\maketitle

\section{Introduction}
\subsection{Homogeneous permutations and the random permutation.} In a paper in 2002, Peter Cameron regarded finite permutations as two linear orders on a finite set, thereby taking a more ``passive'' perspective on permutations than the one which views them as bijections~\cite{CameronPermutations}. He showed that there exist precisely four Fra\"{i}ss\'{e} classes (in the sense of~\cite{Hodges}) of finite permutations in this sense, one of which is the class of \emph{all} finite structures with two linear orders. The Fra\"{i}ss\'{e} limit of the latter class, which is called the \emph{random permutation} and which we denote by $\Pi=(D;<_1,<_2)$, therefore is the (up to isomorphism) unique countable homogeneous structure with two linear orders which contains all finite permutations as induced substructures. Both linear orders of the random permutation are isomorphic to the order of the rational numbers, and the random permutation is the result that appears with probability one in the natural random process that constructs both orders independently. From this it becomes clear that the random permutation cannot correspond to a single bijection on its domain $D$: indeed, it represents a double coset $\Aut(D;<_2)\circ\pi\circ \Aut(D;<_1)$ in the full symmetric group $\Sym(D)$ on $D$, where $\pi$ is any isomorphism from $(D;<_1)$ to $(D;<_2)$, and $\Aut(D;<_i)$ denotes the automorphism group of $(D;<_i)$, for $i=1,2$.

\subsection{Symmetries of the random permutation.} The random permutation possesses two kinds of obvious symmetries. Firstly, it inherits symmetries of the order of the rational numbers: for example, the structure $(D;>_1)$ is obviously isomorphic to $(D;<_1)$, and it is easy to see that likewise $(D;>_1,<_2)$ is isomorphic to $\Pi=(D;<_1,<_2)$. The symmetries of the order of the rational numbers have been classified by Cameron in a famous paper in 1976~\cite{Cameron5}; they are basically composed of two non-trivial symmetries, one of which is reversing the order, and the other one is turning the order cyclically. The second obvious symmetry of $\Pi$ is the fact that not only the orders $(D;<_1)$ and $(D;<_2)$ are isomorphic, but also $(D;<_2,<_1)$ is isomorphic to $\Pi=(D;<_1,<_2)$. 

The symmetries in the above sense of a structure correspond to those subgroups of the full symmetric group of its domain which contain the automorphism group of the structure and which are \emph{closed} in the topology of pointwise convergence. Combining the two kinds of obvious symmetries of $\Pi$ mentioned above, Cameron counted 37 closed supergroups of $\Aut(\Pi)$, and asked whether there were any others, stating the following problem:

\begin{prob}[Problem~2 in~\cite{CameronPermutations}, rephrased]
Determine the closed subgroups of $\Sym(D)$ which contain $\Aut(\Pi)$.
\end{prob}
In this paper, we solve this problem, showing that there exist precisely 39 closed supergroups of $\Aut(\Pi)$. While there turn out be a few groups which had not been considered in~\cite{CameronPermutations}, some of those counted in that paper actually coincide.

\subsection{Reducts and Thomas' conjecture.} For structures $\Gamma, \Delta$ on the same domain, we call $\Gamma$ a  \emph{reduct} of $\Delta$ iff all of its relations and functions have first-order definitions in $\Delta$ without parameters. It follows from the theorem of Ryll-Nardzewski, Engeler, and Svenonius (see e.g.~\cite{Hodges} for all standard model-theoretic notions and theorems) that if we consider two reducts equivalent iff they are reducts of one another, then the reducts of an \emph{$\omega$-categorical} structure $\Delta$ correspond precisely to the closed supergroups of the automorphism group $\Aut(\Delta)$. In this correspondence, every reduct $\Gamma$ of $\Delta$ is sent to $\Aut(\Gamma)$, defining a subjective map onto the closed supergroups of $\Aut(\Delta)$ whose kernel is the above-mentioned equivalence. Since the closed supergroups of $\Aut(\Delta)$ form a complete lattice, so do the reducts of $\Delta$ up to equivalence, the order being provided by first-order definability.

 In 1991, Simon Thomas conjectured that every countable structure which is homogeneous in a finite relational language has only finitely many reducts up to equivalence~\cite{RandomReducts}. At the time, the reducts of only two interesting structures which fall into the scope of the conjecture had been classified: those of the order of the rational numbers (5~reducts)~\cite{Cameron5} and those of the random graph (5~reducts)~\cite{RandomReducts}. Since then the reducts of the random hypergraphs~\cite{Thomas96}, the random tournament~\cite{Bennett-thesis}, the order of the rationals with a constant~\cite{JunkerZiegler}, and more recently those of the random partial order~\cite{Poset-Reducts}, the $K_n$-free graphs with a  constant~\cite{Andras-thesis} and the random ordered graph~\cite{42} have been determined, in all cases confirming Thomas' conjecture. Our classification verifies the conjecture for the random permutation.

\subsection{Superpositions of homogeneous structures.} Let $\C_1, \C_2$ be Fra\"{i}ss\'{e} classes of finite structures in disjoint signatures $\sigma_1$ and $\sigma_2$, respectively, and assume moreover that both classes have \emph{strong amalgamation}. Then the class of finite structures with signature $\sigma_1\cup\sigma_2$ whose restriction to the signature $\sigma_i$ is an element of $\C_i$ for $i=1,2$ is a Fra\"{i}ss\'e class as well. Moreover, the restriction of its Fra\"{i}ss\'e limit $\Delta$ to the signature $\sigma_i$ is the Fra\"{i}ss\'e limit of $\C_i$ for $i=1,2$. In this situation, we say that $\Delta$ is the \emph{free superposition} of the Fra\"{i}ss\'e limits of $\C_1$ and $\C_2$. Using this terminology, the random permutation is the free superposition of two copies of the order of the rational numbers. 

It was only very recently that the reducts of a freely superposed structure, namely the superposition of the random graph and the order of the rational numbers called the \emph{random ordered graph}, were classified up to equivalence~\cite{42}. Our result is the second such classification. One notable contrast between the situation in~\cite{42} and our situation is that the two relations of the random ordered graph are very different, the graph relation being a quite ``free'' binary relation as opposed to the order relation, which gives rise to some asymmetry; in particular, the two relations cannot be flipped. 

In the case of the random permutation, another kind of rather surprising asymmetry appears with respect to possible combinations of the reducts of the two orders. As implied above, one closed supergroup of $\Aut(D;<_1)$ is the one consisting of all order preserving and all order reversing permutations; another one is the one consisting of all permutations which turn the order cyclically. While the first group can be combined with the corresponding group above $\Aut(D;<_2)$ to the group consisting of all permutations which either reverse or preserve both orders simultaneously, the groups of cyclic turns have no similar ``simultaneous'' action -- see the discussion in Section~\ref{section:lattice} for more details.

\subsection{Canonical functions and Ramsey theory.}
We prove our result using a method originally invented in the context of \emph{constraint satisfaction}~\cite{BodPin-Schaefer, RandomMinOps} and further developed in~\cite{BP-reductsRamsey, BPT-decidability-of-definability}. Based on so-called \emph{canonical functions}, this method turned out to be very effective in reduct classifications of homogeneous structures with a \emph{Ramsey expansion}. First applied to this kind of problem in 2011 to determine the reducts of the random partial order~\cite{Poset-Reducts}, it has since served to find the reducts of the $K_n$-free graphs with a  constant~\cite{Andras-thesis} and the 
random ordered graph~\cite{42}. As in the case of the latter structure, we take the approach of first identifying the join irreducible elements of the lattice of closed supergroups of $\Aut(\Pi)$ with the help of canonical functions. We then use canonical functions again to prove that every closed supergroup of $\Aut(\Pi)$ is a join of these groups, exploiting the fact that $\Pi$ is itself a \emph{Ramsey structure} (cf.~Section~\ref{sect:prelims}).

\subsection{A model of the random permutation.} It is helpful to visualize $\Pi$ by means of the following concrete representation of this structure. Let ${\mathbb Q}$ be the rational numbers with the usual order $<$. Call a subset $S$ of ${\mathbb Q}^2$ \emph{independent} iff for all $x,y\in S$ we have $x_1\neq y_1$ and $x_2\neq y_2$. Then the following is easily verified using the fact that $\Pi$ is, up to isomorphism, uniquely determined by the \emph{expansion property}~\cite{Hodges}. 

\begin{fact}\label{fact:model}
Let $D$ be any dense and independent subset of ${\mathbb Q}^2$. Then setting $x<_i y$ iff $x_i<y_i$ for $i=1,2$, we have that $(D;<_1,<_2)$ is a model of (the theory of) $\Pi$.
\end{fact}

\subsection{Acknowledgements.}

The first author would like to thank her advisor \'{A}gnes Szendrei for her continued guidance and support, and for introducing her to the second author and this problem. The second author is indebted to Igor Dolinka and Dragan Ma\v{s}ulovi\'c for drawing his attention to Peter Cameron's question, as well as for valuable discussion and generous hospitality during his visit at the University of Novi Sad. He would also like to thank \'{A}gnes Szendrei  and Keith Kearnes for their equally generous hospitality during his visit at the University of Colorado at Boulder.

\section{The Reducts of $\Pi$}\label{sect:reducts}

\subsection{Generators of closed supergroups of $\Aut(\Pi)$.} With the aim of listing the closed supergroups of $\Aut(\Pi)$, we shall now provide a finite set of  permutations on $D$ such that every closed supergroup of $\Aut(\Pi)$ is \emph{generated} by a subset of that set, in the following sense.

\begin{defn}
Let $\F$ be a set of permutations on $D$, and let $\G$ be a closed permutation group on $D$.
We say that $\F$  \emph{generates} $\G$ (over $\Aut(\Pi)$) iff $\G$ is the smallest closed permutation group that contains $\F\cup \Aut(\Pi)$; in that case, we write $\G=\cl{\F}$. We always assume $\Aut(\Pi)$ to be present in the generating process, and will not mention it explicitly.
When $\F=\{f_1,\ldots,f_n\}$, then we also write $\cl{f_1,\ldots,f_n}$ for $\cl{\F}$.
\end{defn}
The elements of $\cl{\F}$ are precisely those permutations $g$ of $D$ with the property that for all finite $A\subseteq D$ there exists a \emph{term} function over the set $\F\cup \Aut(\Pi)$ which agrees with $g$ on $A$. Here, terms are composites of elements of $\F\cup \Aut(\Pi)$ and of inverses of such elements.

As noted before, the structures $(D;>_1,<_2)$, $(D;<_1,>_2)$, and $(D;>_1,>_2)$ are all isomorphic to $\Pi$. Let $\mix \rev\id$, $\mix \id\rev$ and $\mix \rev\rev$ be isomorphisms from $\Pi$ to these structures: that is, $\mix \rev\id$ reverses $<_1$ while preserving $<_2$, $\mix \id\rev$ does the same with the roles of the two orders interchanged, and $\mix \rev\rev$ reverses both orders. Moreover, $(D;<_2,<_1)$ is isomorphic to $\Pi$; let $\sw$ be an isomorphism. 

In the model of $\Pi$ provided in Fact~\ref{fact:model}, we can visualize these permutations as follows. Observe that if $D'\subseteq {\mathbb Q}^2$ is dense and independent, then there exist automorphisms $\alpha_1,\alpha_2$ of $({\mathbb Q};<)$ such that $\alpha:=(\alpha_1,\alpha_2):{\mathbb Q}^2\To{\mathbb Q}^2$ maps $D'$ bijectively onto $D$. Moreover, if automorphisms $\beta_1,\beta_2$ of $({\mathbb Q};<)$ are so that $\beta:=(\beta_1,\beta_2):{\mathbb Q}^2\To{\mathbb Q}^2$ maps $D'$ bijectively onto $D$, then there exists $\gamma\in\Aut(\Pi)$ such that $\alpha=\gamma\circ\beta$. Hence, every function $f: {\mathbb Q}^2\To {\mathbb Q}^2$ with the property that it sends $D$ bijectively onto a set $D'$ which is dense and independent induces permutations on $D$ of the form $\alpha\circ f\rest_D$, and any two permutations of this form are equivalent for our purposes since they generate the same closed groups.

In this construction, $\mix \rev\id$ is induced by the mapping from ${\mathbb Q}^2$ to ${\mathbb Q}^2$ which sends any $(x_1,x_2)$ to $(-x_1,x_2)$; we may thus say that geometrically, $\mix \rev\id$ corresponds to the mapping $(x_1,x_2)\mapsto (-x_1,x_2)$ on ${\mathbb Q}^2$. Similarly, $\mix \id\rev$ corresponds to $(x_1,x_2)\mapsto (x_1,-x_2)$, and $\mix \rev\rev$ to $(x_1,x_2)\mapsto (-x_1,-x_2)$, which is just the composite of the preceding two functions. The function $\sw$ is geometrically nothing else but $(x_1,x_2)\mapsto (x_2,x_1)$.

We use our model of $\Pi$ in order to define more permutations. Let $r\in {\mathbb R}\setminus {\mathbb Q}$ be an irrational number, and let $f_r$ be any function which sends the interval $(-\infty,r)\cap{\mathbb Q}$ bijectively onto $(r,\infty)\cap{\mathbb Q}$ whilst preserving the order on $(-\infty,r)\cap \Q$ and $(r,\infty)\cap\Q$. Then $(x_1,x_2)\mapsto (f_r(x_1),x_2)$ is a permutation of ${\mathbb Q}^2$ which induces a permutation on $\Pi$ as described above -- we denote this permutation by $\mix {t_r} \id$. It is straightforward to see that the closed group generated by such a function is independent of $r$, and we will thus write $\mix {t} \id$ whenever there is no need to refer to $r$ explicitly. Similarly, we define functions $\mix \id {t_r}$ and $\mix \id{t}$.

\subsection{Closed supergroups of $\Aut(D;<_i)$.} Recall that $(D;<_i)$ is isomorphic to the order of the rational numbers, and that the closed supergroups of the automorphism group of that order have been classified~\cite{Cameron5}. In our context, that classification can be stated as follows.

\begin{thm}[Cameron~\cite{Cameron5}]\label{thm:cameron5}
The closed supergroups of $\Aut(D;<_1)$ are precisely the following:
\begin{enumerate}
\item $\Aut(D;<_1)$;
\item $\cl{\{\mix \rev\id \}\cup \Aut(D;<_1)}$;
\item $\cl{\{\mix t\id \}\cup \Aut(D;<_1)}$;
\item $\cl{\{\mix \rev\id, \mix t\id \}\cup \Aut(D;<_1)}$;
\item $\Sym(D)$.
\end{enumerate}
\end{thm}

Of course, the theorem for $(D;<_2)$ is similar. If we wish to see these groups as automorphism groups of reducts of $(D;<_i)$, then the following relations on $D$ are suitable. For $i\in\{1,2\}$, set
\begin{itemize}
\item $\btw_i(x,y,z)\Leftrightarrow(x<y<z)\vee(z<y<x)$;
\item $\cyc_i(x,y,z)\Leftrightarrow(x<y<z)\vee(y<z<x)\vee(z<x<y)$;
\item $\sepa_i(w,x,y,z)\Leftrightarrow(\cyc_i(w,x,y)\wedge\cyc_i(w,z,x))\vee(\cyc_i(w,y,x)\wedge\cyc_i(w,x,z))$.
\end{itemize}

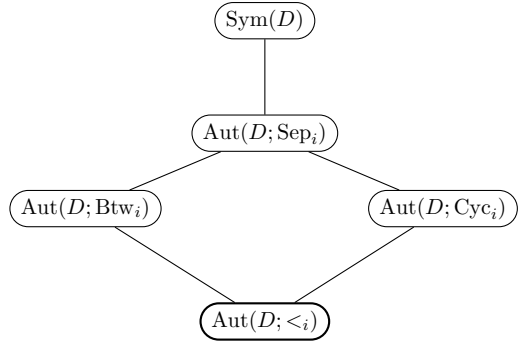
\begin{figure}
\begin{tikzpicture}[every node/.style={fill=white}]

	\node(i) [JI] {$\Aut(D;<_i)$};
	
	\node(ci) [subgroup, above left=of i]{$\Aut(D;\btw_i)$};
	\path[line] (i) -- (ci);
	
	\node(di) [subgroup, above right=of i]{$\Aut(D;\cyc_i)$};
	\path[line] (i) -- (di);
	
	\node(cdi) [subgroup, above=2cmof i]{$\Aut(D;\sepa_i)$};
	\path[line] (ci) -- (cdi);
	\path[line] (di) -- (cdi);
	
	\node(top) [subgroup, above=of cdi]{$\Sym(D)$};
	\path[line] (cdi) -- (top);
	
\end{tikzpicture}
\caption{Closed supergroups of $\Aut(D;<_i)$.}
\end{figure}

\begin{cor}[Cameron~\cite{Cameron5}] \label{cor:above aut}
The closed supergroups of $\Aut(D;<_1)$ are precisely the following:
\begin{enumerate}
\item $\Aut(D;<_1)$;
\item $\Aut(D;\btw_1)$;
\item $\Aut(D;\cyc_1)$;
\item $\Aut(D;\sepa_1)$.
\item $\Aut(D;=)$.
\end{enumerate}
\end{cor}

The groups in Theorem~\ref{thm:cameron5} and Corollary~\ref{cor:above aut} are listed in the same order.

\subsection{Join irreducible closed supergroups of $\Aut(\Pi)$.} Arbitrary intersections of closed permutations groups on $D$ yield closed permutation groups. Therefore, the closed permutation groups on $D$ form a complete lattice with respect to inclusion, and the closed supergroups of $\Aut(\Pi)$ form an interval $\mathfrak L$ therein. We now provide the set of all completely join irreducible elements of the lattice $\mathfrak L$, i.e., of all elements of $\mathfrak L$ which are not the (in theory, possibly infinite) join of other groups in $\mathfrak L$.

\begin{defn} 
Let $\JI$ consist of the following groups:
\begin{enumerate}
\item $\langle\mix\id\rev \rangle$;
\item $\langle\mix\id {t} \rangle$;
\item $\langle\mix\rev\id \rangle$;
\item $\langle\mix{t} \id \rangle$;
\item $\langle\mix \rev\rev\rangle$;
\item $\langle\sw\rangle$;
\item $\langle\sw\circ\mix\rev\rev \rangle$;
\item $\langle\sw\circ\mix\id\rev \rangle$;
\item $\Aut(D;\oo)$;
\item $\Aut(D;\ot)$.
\end{enumerate}
\end{defn}

We are going to prove the following theorem, which implies that the closed permutation groups which properly contain $\Aut(\Pi)$ are precisely the joins of groups in $\JI$. As a consequence, it follows that there are at most $2^{|\JI|}+1=2^{10}+1$ closed supergroups of $\Aut(\Pi)$.

\begin{thm}\label{thm:JI}
Let $\G\supseteq\Aut(\Pi)$ be a closed group and let $f\in\Sym(D)$ be such that $f\notin \G$. Then there exists a group $\H\in \JI$ such that $\H\subseteq\langle\{f\}\cup \G\rangle$ and $\H\nsubseteq \G$.
\end{thm}

\begin{cor}\label{cor:JI}
Let $\G\supsetneq\Aut(\Pi)$ be a closed group. Then $\G$ is the join of elements of $\JI$. In particular, $\mathfrak L$ is finite.
\end{cor}

By systematically investigating the joins of elements of $\JI$, we then obtain that there exist precisely 39 distinct closed supergroups of $\Aut(\Pi)$, and determine the exact shape of $\mathfrak L$. In order to show a compact picture of $\mathfrak L$, we name the elements of $\JI$ as follows. First those which we know from the classification of the symmetries of the order of the rational numbers\ldots

\begin{table}[h]
\begin{tabular}{|c|c|c|c|c|c}\hline
Letter & a & b & c & d & e  \\\hline
Group  & $\la\mix\id\rev\ra$ & $\la\mix\id{t}\ra$ & $\la\mix\rev\id\ra$ & $\la\mix{t}\id\ra$ & $\la\mix\rev\rev\ra$\\\hline
\end{tabular}
\end{table}

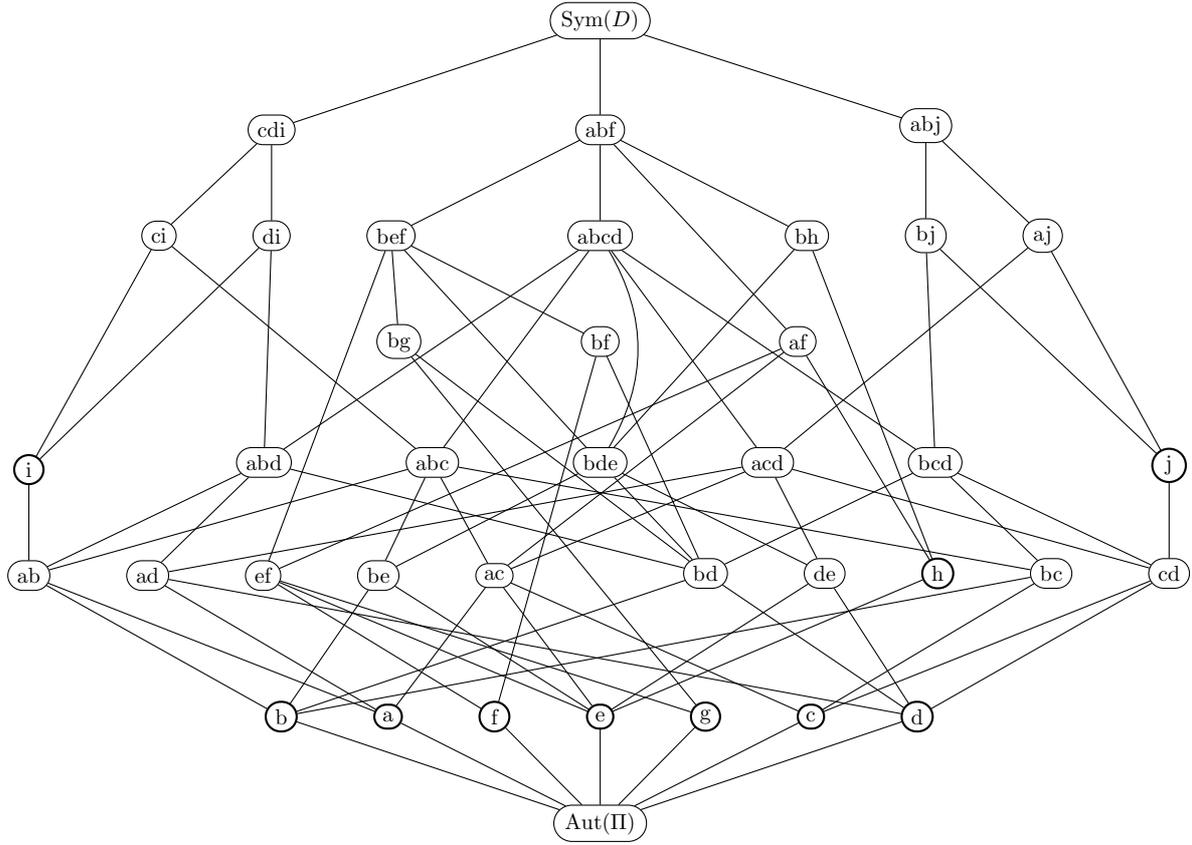
\begin{figure}
\begin{tikzpicture}[every node/.style={fill=white}]

\node(0)[subgroup, ]{$\Aut(\Pi)$};

\node(e)[JI,above=of 0 ]{e};
\path[line](0)--(e);

\node(f)[JI,left=of e]{f};
\path[line](0)--(f);

\node(a)[JI,left=of f]{a};
\path[line](0)--(a);

\node(b)[JI,left=of a]{b};
\path[line](0)--(b);

\node(g)[JI,right=of e]{g};
\path[line](0)--(g);

\node(c)[JI,right=of g]{c};
\path[line](0)--(c);

\node(d)[JI,right=of c]{d};
\path[line](0)--(d);

\node(ac)[subgroup,above=1.5cmof f]{ac};
\path[line](a)--(ac);
\path[line](c)--(ac);
\path[line](e)--(ac);

\node(be)[subgroup, left=of ac]{be};
\path[line](b)--(be);
\path[line](e)--(be);

\node(ef)[subgroup,left=of be]{ef};
\path[line](e)--(ef);
\path[line](f)--(ef);
\path[line](g)--(ef);

\node(ad)[subgroup,left=of ef]{ad};
\path[line](a)--(ad);
\path[line](d)--(ad);

\node(ab)[subgroup,left=of ad]{ab};
\path[line](a)--(ab);
\path[line](b)--(ab);

\node(bd)[subgroup,above=1.5cmof g]{bd};
\path[line](b)--(bd);
\path[line](d)--(bd);

\node(de)[subgroup,right=of bd]{de};
\path[line](d)--(de);
\path[line](e)--(de);

\node(h)[JI,right=of de]{h};
\path[line](e)--(h);

\node(bc)[subgroup,right=of h]{bc};
\path[line](b)--(bc);
\path[line](c)--(bc);

\node(cd)[subgroup,right=of bc]{cd};
\path[line](c)--(cd);
\path[line](d)--(cd);

\node(bde)[subgroup,above=3cmof e]{bde};
\path[line](bd)--(bde);
\path[line](be)--(bde);
\path[line](de)--(bde);

\node(abc)[subgroup,left=1.5cmof bde]{abc};
\path[line](ab)--(abc);
\path[line](ac)--(abc);
\path[line](bc)--(abc);
\path[line](be)--(abc);

\node(abd)[subgroup,left=1.5cmof abc]{abd};
\path[line](ab)--(abd);
\path[line](ad)--(abd);
\path[line](bd)--(abd);

\node(acd)[subgroup,right=1.5cmof bde]{acd};
\path[line](ac)--(acd);
\path[line](ad)--(acd);
\path[line](cd)--(acd);
\path[line](de)--(acd);

\node(bcd)[subgroup,right=1.5cmof acd]{bcd};
\path[line](bc)--(bcd);
\path[line](bd)--(bcd);
\path[line](cd)--(bcd);

\node(bf)[subgroup,above=1.2cmof bde]{bf};
\path[line](f)--(bf);
\path[line](bd)--(bf);

\node(bg)[subgroup,left=2.1cmof bf]{bg};
\path[line](g)--(bg);
\path[line](bd)--(bg);

\node(af)[subgroup,right=2.1cmof bf]{af};
\path[line](h)--(af);
\path[line](ac)--(af);
\path[line](ef)--(af);

\node(i)[JI,above=of ab]{i};
\path[line](ab)--(i);

\node(j)[JI,above=of cd]{j};
\path[line](cd)--(j);

\node(abcd)[subgroup, above=of bf]{abcd};
\path[line](abc)--(abcd);
\path[line](abd)--(abcd);
\path[line](acd)--(abcd);
\path[line](bcd)--(abcd);
\path[line](bde) edge[bend right] (abcd);

\node(bef)[subgroup,left=2cmof abcd]{bef};
\path[line](bf)--(bef);
\path[line](bg)--(bef);
\path[line](ef)--(bef);
\path[line](bde)--(bef);

\node(bh)[subgroup,right=2cmof abcd]{bh};
\path[line](h)--(bh);
\path[line](bde)--(bh);

\node(di)[subgroup,left=of bef]{di};
\path[line](i)--(di);
\path[line](abd)--(di);

\node(bj)[subgroup,right=of bh]{bj};
\path[line](j)--(bj);
\path[line](bcd)--(bj);

\node(aj)[subgroup,right=of bj]{aj};
\path[line](j)--(aj);
\path[line](acd)--(aj);

\node(ci)[subgroup,left=of di]{ci};
\path[line](i)--(ci);
\path[line](abc)--(ci);

\node(abf)[subgroup,above=of abcd]{abf};
\path[line](af)--(abf);
\path[line](bh)--(abf);
\path[line](bef)--(abf);
\path[line](abcd)--(abf);

\node(abj)[subgroup,above=of bj]{abj};
\path[line](aj)--(abj);
\path[line](bj)--(abj);

\node(cdi)[subgroup,above=of di]{cdi};
\path[line](ci)--(cdi);
\path[line](di)--(cdi);

\node(1)[subgroup,above=of abf]{$\Sym(D)$};
\path[line](abf)--(1);
\path[line](abj)--(1);
\path[line](cdi)--(1);

\end{tikzpicture}
\caption{The lattice $\mathfrak L$ of closed supergroups of $\Aut(\Pi)$.}\label{fig:lattice}
\end{figure}

\noindent \ldots and then those which we get by switching the orders, or by completely ignoring one of the orders. In Figure~\ref{fig:lattice}, each group in $\mathfrak L$ is labeled by a minimal set of elements of $\JI$ whose join it equals.

\begin{table}[h]
\begin{tabular}{c|c|c|c|c|}\hline
 f & g & h & i & j \\\hline
  $\la\sw\ra$ & $\la\sw\circ\mix\rev\rev\ra$ & $\la\sw\circ\mix\id\rev\ra$ & $\Aut(\oo)$ &  $\Aut(\ot)$\\\hline
\end{tabular}
\end{table}

\begin{thm}\label{thm:L}
The lattice $\mathfrak L$ of closed supergroups of $\Aut(\Pi)$ has 39~elements, and the shape as represented in Figure~\ref{fig:lattice}.
\end{thm}

\subsection{Organization of the paper.} In the following section (Section~\ref{sect:prelims}), we provide the Ramsey-theoretic preliminaries for the proof of Theorem~\ref{thm:JI}, consisting of a Ramsey-type statement for the class of finite permutations, and the method of \emph{canonical functions}. After that, in Section~\ref{section:proof}, we give the proof of Theorem~\ref{thm:JI}. In the final section (Section~\ref{section:lattice}), we show which joins of element in $\JI$ actually coincide, from which we obtain the precise size and shape of $\L$, proving Theorem~\ref{thm:L}. It is also there that we discuss the shape of $\L$ in more detail, and compare it with Peter Cameron's count from~\cite{CameronPermutations}.

\section{Ramsey-theoretic Preliminaries.}\label{sect:prelims}

\subsection{Ramsey structures.} Our proof exploits a Ramsey-type property of $\Pi$, as in the following definition.

\begin{defn}\label{def:ramseystructure}
Let $\Delta$ be a countable relational structure with domain $X$. 
For any structure $\Gamma$ in the language of $\Delta$, write $\mix \Delta \Gamma$ for the set of all induced substructures of $\Delta$ which are isomorphic to $\Gamma$. Then $\Delta$ is called a \emph{Ramsey structure} iff for all finite induced substructures $\Omega$ of $\Delta$, all induced substructures $\Gamma$ of $\Omega$, and all $\chi\colon\mix \Delta \Gamma\To 2$ there exists $\Omega'\in \mix \Delta \Omega$ such that the restriction of $\chi$ to $\mix {\Omega'} \Gamma$ is constant.
\end{defn}

For example, the order of the rational numbers is Ramsey: this fact is easily seen to be equivalent to Ramsey's theorem. It follows that $\Pi$ is Ramsey, since it is the free superposition of two copies of a homogeneous relational Ramsey structure whose finite substructures are \emph{rigid} and have strong amalgamation. That such superpositions are Ramsey has been proven in~\cite{Bod-New-Ramsey-classes} using infinitary methods, and later in~\cite{Sok-Directed-graphs-and} using finite combinatorics. The fact that $\Pi$ is Ramsey was, however, proven before in~\cite{BoettcherFoniok} and independently in~\cite{Sokic}.

\begin{fact}
$\Pi$ is a Ramsey structure.
\end{fact}

\subsection{Canonical functions.}

The fact that $\Pi$ is a relational homogeneous Ramsey structure implies that distinct closed supergroups of $\Aut(\Pi)$ can be distinguished by  so-called \emph{canonical functions}. This has been observed in~\cite{BP-reductsRamsey, BPT-decidability-of-definability}, and will be our method for proving our main result.

\begin{defn} 
Let $\Delta$ be a structure, and let $a$ be an $n$-tuple of elements in $\Delta$. The \emph{type} of $a$ in $\Delta$ is the set of first-order formulas with free variables $x_1,\ldots,x_n$ that hold for $a$ in $\Delta$.
\end{defn}

\begin{defn} 
Let $\Delta$ and $\Gamma$ be structures. A \emph{type condition} between $\Delta$ and $\Gamma$ is a pair $(t,s)$, such that $t$ is the type on an $n$-tuple in $\Delta$ and $s$ is the type of an $n$-tuple in $\Gamma$, for some $n\geq 1$. A function $f:\Delta\to\Gamma$ \emph{satisfies} a type condition $(t,s)$ iff the type of $f(a)$ in $\Gamma$ equals $s$ for all $n$-tuples $a$ in $\Delta$ of type $t$. 

A \emph{behavior} $\mathcal{B}$ is a set of type conditions between $\Delta$ and $\Gamma$, and a function $f:\Delta\to\Gamma$ \emph{has behavior} $\mathcal{B}$ iff it satisfies all type conditions in $\mathcal{B}$.
\end{defn}

\begin{defn} 
Let $\Delta$ and $\Gamma$ be structures. A function $f:\Delta\to\Gamma$ is \emph{canonical} iff for every type $t$ of an $n$-tuple in $\Delta$ there is a type $s$ of an $n$-tuple in $\Gamma$ such that $f$ satisfies the type condition $(t,s)$. That is, canonical functions send $n$-tuples of the same type to $n$-tuples of the same type, for all $n\geq 1$.
\end{defn}

Note that any canonical function induces a function from the types over $\Delta$ to the types over $\Gamma$. The canonical functions that are needed are in general not permutations: they fail to be surjective. In fact, we obtain our classification of the closed supergroups of $\Aut(\Pi)$ by an analysis of the closed transformation monoids of injective functions on $D$ which contain $\Aut(\Pi)$. Here, just as for permutation groups, ``closed'' means closed in the topology of \emph{pointwise convergence}, i.e., the topology of the product space $D^D$, where $D$ is taken to be discrete. The fact that we have to leave the realm of permutation groups necessitates the following definition.

\begin{defn} Let $\mathcal{F}\subseteq D^D$. We say that $\mathcal{F}$ \emph{mon-generates} a function $g\colon D\to D$ (over $\Aut(\Pi)$) iff $g$ is contained in the smallest closed submonoid of $D^D$ which contains $\F\cup\Aut(\Pi)$. In other words, this is the case iff
for every finite subset $A\subseteq D$ there exist $n\geq 1$ and $f_1,\ldots, f_n\in\mathcal{F}\cup\Aut(\Pi)$ such that $f_1\circ\cdots\circ f_n$ agrees with $g$ on $A$. In this paper, when we write that $\mathcal{F}$ mon-generates a function then we always mean ``over $\Aut(\Pi)$''. For functions $f,g\colon D\to D$ we shall say that $f$ mon-generates $g$ rather than $\{f\}$ mon-generates $g$.
\end{defn}

Our proof relies on the following proposition which is a consequence of~\cite{BP-reductsRamsey, BPT-decidability-of-definability} and the fact that $\Pi$ is a Ramsey structure. For $c_1,\ldots,c_n\in D$, let $(\Pi,c_1,\ldots,c_n)$ denote the structure obtained from $\Pi$ by adding the constants $c_1,\ldots,c_n$ to the language.

\begin{prop}\label{prop:magical}
Let $f:D\To D$ be any injective function, and let $c_1,\ldots,c_n\in D$. Then $f$ mon-generates an injective function $g:D\To D$ such that
\begin{itemize}
\item $g$ agrees with $f$ on $\{c_1,\ldots,c_n\}$;
\item $g$ is canonical as a function from $(\Pi,c_1,\ldots,c_n)$ to $\Pi$.
\end{itemize}
\end{prop}

Note that any two canonical functions from $(\Pi,c_1,\ldots,c_n)$ to $\Pi$ with equal behavior, i.e., which satisfy the same type conditions, mon-generate one another: this is a consequence of the homogeneity of $\Pi$. Thus, for a fixed choice  of $c_1,\ldots,c_n$, there are essentially only finitely many distinct canonical functions, and they are essentially finite objects since they are determined by their behavior. 

Our proof of Theorem~\ref{thm:JI} uses the following idea: if $\G\in\JI$ and $f\in\Sym(D)\setminus \G$, then there exist  $c_1,\ldots,c_n\in D$ such that no function in $\G$ agrees with $f$ on $\{c_1,\ldots,c_n\}$. Let $g$ be the canonical function mon-generated by $f$ by virtue of Proposition~\ref{prop:magical}; then $g$ is not mon-generated by $\G$. By analyzing the possible behaviors of $g$, we deduce that $g$ mon-generates all functions of some $\H\in\JI$ which is not contained in $\G$. But then $\cl{\{f\}\cup \G}$ contains $\H$, and hence a new element of $\JI$.

\section{The Proof}\label{section:proof}

\subsection{All canonical functions from $\Pi$ to $\Pi$.}\label{subsection:canonical} We start by investigating all behaviors of canonical injections from $\Pi$ to $\Pi$. As a matter of fact, this will make us rediscover many of the functions presented in Section~\ref{sect:reducts}. 

\begin{defn}
	Let $g,h\colon\Pi\To\Pi$ be canonical. We say that $g$ \emph{behaves like $h$} iff for all $x,x'\in X$ we have that the type of $(g(x),g(x'))$ equals the type of $(h(x),h(x'))$ in $\Pi$.
\end{defn}

\begin{defn}
Define the following binary relations on $D$:
\begin{itemize}
\item $\up(x,y)\Leftrightarrow x\oo y\wedge x\ot y$;
\item $\st(x,y)\Leftrightarrow\up(x,y)\vee \up(y,x)$;
\item $\tw(x,y)\Leftrightarrow\neg\st(x,y)$;
\item $\dow(x,y)\Leftrightarrow x\oo y\wedge y\ot x$.
\end{itemize}
We call a subset $S\subseteq D$ \emph{diagonal} iff either $\st(x,y)$ for all distinct $x,y\in S$, or $\tw(x,y)$ for all distinct $x,y\in S$.
\end{defn}

\begin{prop}\label{prop:canonical} 
Let $\G$ be a closed supergroup of $\Aut(\Pi)$ and let $g\colon \Pi\To\Pi$ be a canonical function mon-generated by $\G$. Then either the image of $g$ is diagonal and $\Aut(D;<_i)\subseteq \G$ for some $i\in\{1,2\}$, or $g$ behaves like one of the following functions:
\begin{enumerate}
\item the identity function $\id$ on $D$;
\item $\mix\id\rev$;
\item $\mix\rev\id$;
\item $\mix\rev\rev$;
\item $\sw$;
\item $\sw\circ\mix\rev\rev$;
\item $\sw\circ\mix\id\rev$;
\item $\sw\circ\mix\rev\id$.
\end{enumerate}
\end{prop}
\begin{proof} 
There are precisely four types of pairs of distinct elements in $\Pi$. Let $t_1,t_2,t_3$ and $t_4$ be the types of a pair $(x,y)$ with $\up(x,y)$, $\dow(x,y)$, $\up(y,x)$, and $\dow(y,x)$, respectively. The behavior of $g$ is fully specified by how it behaves on pairs of types $t_1$ and $t_2$. This gives 16 possible canonical behaviors. 

Take any $x,y,u,v\in D$ such that the types of $(x,y)$ and $(u,v)$ in $\Pi$ are $t_1$ and $t_2$, respectively. Suppose first that $(g(x),g(y))$ and $(g(u),g(v))$ both have type $t_1$. Then $g$ sends $D$ to a diagonal set while preserving $\oo$. Take any $a_1,\ldots,a_n,b_1,\ldots,b_n\in D$ with $a_i\oo a_{i+1}$ and $b_i\oo b_{i+1}$ for all $1\leq i\leq n-1$. In order to show that $\Aut(D;\oo)\subseteq \G$, it suffices to show that there is a function $f$ mon-generated by $\G$ such that $f(a_i)=b_i$ for all $1\leq i\leq n$.

Since $g$ is canonical, for all $1\leq i\leq n-1$, the pairs $(g(a_i),g(a_{i+1}))$ and $(g(b_i),g(b_{i+1}))$ have type $t_1$. Therefore, the tuples $(g(a_1),\ldots,g(a_n))$ and $(g(b_1),\ldots,g(b_n))$ have the same type in $\Pi$. Since $\Pi$ is $\omega$-categorical, there exists $\alpha\in\Aut(\Pi)$ such that $\alpha(g(a_i))=g(b_i)$ for all $1\leq i\leq n$. Moreover, since $g$ is mon-generated by $\G$, there exists $h\in\G$ which agrees with $g$ on $\{b_1,\ldots,b_n\}$. Let $f=h^{-1}\circ\alpha\circ g$. Then $f$ is a function mon-generated by $\G$ such that $f(a_i)=b_i$ for all $1\leq i\leq n$, giving that $\Aut(D;\oo)\subseteq \G$.

If $(g(x),g(y))$ and $(g(w),g(z))$ both have type $t_2$, then $g$ also sends $D$ to a diagonal set while preserving $\oo$. Similarly, if $(g(x),g(y))$ and $(g(w),g(z))$ both have type $t_3$ or $t_4$, then $g$ sends $D$ to a diagonal set while reversing $\oo$. By the same argument as above, $\Aut(D;\oo)\subseteq \G$ holds in these cases as well.

Now suppose $(g(x),g(y))$ has type $t_1$ and $(g(w),g(z))$ has type $t_3$. Then $g$ sends $D$ to a diagonal set while preserving $\ot$. This is also the case if $(g(x),g(y))$ has type $t_4$ and $(g(w),g(z))$ has type $t_2$. Similarly, if $(g(x),g(y))$ has type $t_3$ and $(g(w),g(z))$ has type $t_1$, or if $(g(x),g(y))$ has type $t_2$ and $(g(w),g(z))$ has type $t_4$, then $g$ sends $D$ to a diagonal set while reversing $\ot$. Using a similar argument as above, in each of these cases we see that $\Aut(D;\ot)\subseteq \G$. 

This leaves 8 remaining canonical behaviors, corresponding to the behaviors of the 8 functions listed above.
\end{proof}

\subsection{Canonical behaviors with constants.}\label{subsection:constants}

\begin{defn} 
Let $c_1,\ldots,c_n\in D$ be distinct. For $0\leq i\leq n$, define the \emph{$i$-th column} of $(\Pi,c_1,\ldots,c_n)$, denoted $C_i$, to be the set of all $d\in D\setminus\{c_1,\ldots, c_n\}$ such that exactly $i$ of $\{c_1,\ldots,c_n\}$ are less than $d$ with respect to the order $\oo$. Similarly, define the \emph{$i$-th row} of $(\Pi,c_1,\ldots,c_n)$, denoted $R_i$, to be the set of all $d\in D\setminus\{c_1,\ldots, c_n\}$ such that exactly $i$ of $\{c_1,\ldots,c_n\}$ are less than $d$ with respect to the order $\ot$.
\end{defn}
 
For distinct $c_1,\ldots,c_n\in D$, the infinite orbits of $\Aut(\Pi,c_1,\ldots,c_n)$ are precisely the sets $R_i\cap C_j,$ for all $0\leq i,j\leq n$. From our model of $\Pi$ in Fact~\ref{fact:model}, the following becomes evident.

\begin{fact} Each infinite orbit of $\Aut(\Pi,c_1,\ldots,c_n)$ is isomorphic to $\Pi$.
\end{fact}

By Proposition~\ref{prop:magical}, we know that closed groups containing $\Aut(\Pi)$ can be distinguished by functions which are canonical from $(\Pi,c_1,\ldots,c_n)$ to $\Pi$, for some $c_1,\ldots,c_n\in D$. Every such function behaves like a canonical function from $\Pi$ to $\Pi$ on each infinite orbit of $\Aut(\Pi,c_1,\ldots,c_n)$, in the following sense.

\begin{defn}
	Let $h\colon\Pi\To\Pi$ be canonical, let $X,Y\subseteq D$ be disjoint, and let $g\colon  D\To D$ be a function. We say that $g$ \emph{behaves like $h$ on $X$} iff for all $x,x'\in X$ we have that the type of $(g(x),g(x'))$ equals the type of $(h(x),h(x'))$ in $\Pi$. We say that $g$ \emph{behaves like $h$ between $X,Y$} iff for all $x\in X$, $y\in Y$ we have that the type of $(g(x),g(y))$ equals the type of $(h(x),h(y))$ in $\Pi$.
\end{defn}

In this subsection we will show that in a sense, there are only few relevant canonical functions from $(\Pi,c_1,\ldots,c_n)$ to $\Pi$. More precisely, we will prove Lemma~\ref{lem:permute} which states that 
if $g$ is a canonical injection from $(\Pi,c_1,\ldots,c_n)$ to $\Pi$ which behaves like $\id$ on some infinite orbit of $\Aut(\Pi,c_1,\ldots,c_n)$, then either $g$ is mon-generated by $\la\mix\id{t},\mix{t}\id\ra$, or else any closed group which mon-generates $g$ contains $\Aut(D;<_i)$ for some $i\in\{1,2\}$. From this the proof of Theorem~\ref{thm:JI} quickly follows.

\begin{defn}
	Let $R\subseteq D^k$ be a relation, and let $X_1,\ldots,X_k\subseteq D$. We say that \emph{$R(X_1,\ldots,X_k)$ holds} iff $R(x_1,\ldots,x_k)$ holds for all $x_i\in X_i$.
\end{defn}	
 
\begin{defn}
Let $X,Y\subseteq D$ be disjoint. We say that $X,Y$ are in \emph{diagonal position} iff either $\st(X,Y)$ or $\tw(X,Y)$ holds. For a function $g\colon D\To D$, we say that $g$ \emph{diagonalizes} $X,Y$ iff $g[X],g[Y]$ are in diagonal position.
\end{defn}

\begin{defn}
	Let $X, Y$ be subsets of $D$, $i\in\{1,2\}$, and let $g\colon D\To D$ be a function. We say that $g$ 
\begin{enumerate}
	\item \emph{preserves $<_i$ on $X$} iff $x_1<_ix_2$ implies $g(x_1)<_i g(x_2)$, for all $x_1,x_2\in X$;
	\item \emph{reverses $<_i$ on $X$} iff $x_1<_ix_2$ implies $g(x_2)<_i g(x_1)$, for all $x_1,x_2\in X$;
	\item \emph{preserves $<_i$ between $X$ and $Y$} iff $x<_i y$ implies $g(x)<_i g(y)$, for all $x\in X, y\in Y$; 
	\item \emph{reverses $<_i$ between $X$ and $Y$} iff $x<_i y$ implies $g(y)<_i g(x)$, for all $x\in X, y\in Y$.
\end{enumerate}
\end{defn}

\begin{lem}\label{lem:orbits} 
Let $c_1,\ldots, c_n\in D$ and let $g\colon (\Pi,c_1,\ldots, c_n)\To\Pi$ be canonical. Let $i,j\in\{1,2\}$ with $i\neq j$, and let $X,Y$ be infinite orbits of $\Aut(\Pi, c_1\ldots, c_n)$ with $X<_i Y$. If $g$ behaves like $\id$ on $X$ or on $Y$, then
\begin{enumerate}
\item $g$ either preserves or reverses $<_i$ between $X$ and $Y$;
\item $g$ either preserves $<_j$ on $X\cup Y$ or $g$ diagonalizes $X$ and $Y$.
\end{enumerate}
\end{lem}
\begin{proof} 
Without loss of generality, suppose $i=1$ and that $g$ behaves like $\id$ on $X$. If $X$ and $Y$ are in diagonal position, then all pairs $(x,y)$ with $x\in X$ and $y\in Y$ have the same type in $\Pi$. Therefore, all pairs $(x',y')$ with $x'\in g[X]$ and $y'\in g[Y]$ have the same type in $\Pi$. It follows that $g[X]$ and $g[Y]$ are in diagonal position, verifying~(2), and moreover~(1) holds. We may thus henceforth assume that $X$ and $Y$ are not in diagonal position.

We first show~(1). If $g$ does not preserve $\oo$ between $X$ and $Y$, then there exist $x\in X$ and $y\in Y$ such that $g(y)\oo g(x)$; without loss of generality we may assume $\up(x,y)$. Pick $z\in X$ such that $\dow(z,y)$ and $x\oo z$. Since $g$ preserves $\oo$ on $X$, $g(x)\oo g(z)$. So, by transitivity, $g(y)\oo g(z)$. Since $g$ is canonical it follows that $g$ reverses $\oo $ between $X$ and $Y$. Hence, $g$ either preserves or reverses $\oo$ between $X$ and $Y$.

To see~(2), suppose first that $g$ violates $\ot$ between $X$ and $Y$. Without loss of generality, there exist $x\in X, y\in Y$ with $x\ot  y$ and $g(y)\ot  g(x)$. Take any $u\in X, v\in Y$. If $u\ot v$, then since $g$ is canonical, $g(v)\ot g(u)$. If $v\ot u$, pick $w\in X$ with $w\ot v$. Then by the above, $g(v)\ot  g(w)$. Moreover, since $g$ preserves $\ot$ on $X$, $g(w)\ot g(u)$. So, by transitivity, $g(v)\ot g(u)$. Hence, $g[Y]\ot g[X]$. This together with (1) shows that $g[X]$ and $g[Y]$ are in diagonal position.

Now assume that $g$ preserves $\ot$ between $X$ and $Y$. Take any $x,x'\in X$, $y,y'\in Y$ such that $x\ot y\ot x'\ot y'$. Since $g$ preserves $\ot$ between $X$ and $Y$, $g(x)\ot g(y)\ot g(x')\ot g(y')$. Hence, $g$ preserves $\ot$ on $X\cup Y$.
\end{proof}

\begin{lem}\label{lem:diagonalize}
Let $\G$ be a closed supergroup of $\Aut(\Pi)$. Let $c_1,\ldots, c_n\in D$ and let $g$ be a function mon-generated by $\G$ which is canonical as a function from $(\Pi, c_1,\ldots, c_n)$ to $\Pi$. Suppose $g$ diagonalizes infinite orbits $X,Y$ of $\Aut(\Pi,c_1,\ldots, c_n)$ which are not in diagonal position. Then $\Aut(D;<_i)\subseteq \G$, for some $i\in\{1,2\}$.
\end{lem}
\begin{proof} 
Without loss of generality, suppose $X\oo Y$ and $g[X]\ot g[Y]$. We claim that $g$ mon-generates a function $g'$ which behaves like $\id$ on $X$ and $Y$, and such that $\up(g'[X],g'[Y])$ holds.

Since $g$ is canonical as a function $(\Pi,c_1,\ldots,c_n)\To\Pi$, it behaves like a canonical function from $\Pi$ to $\Pi$ on each of $X$ and $Y$. Say $g$ behaves like $h\colon \Pi\To\Pi$ on $X$ and $k\colon \Pi\To\Pi$ on $Y$. Then $h$ (and similarly $k$) is mon-generated by $\G$: any self-embedding $\iota$ of $\Pi$ whose range is contained in $X$ is mon-generated by $\Aut(\Pi)$, and $g\circ \iota$ behaves like $h$.  
Therefore, by Proposition~\ref{prop:canonical}, either $\Aut(D;<_i)\subseteq \G$ for some $i\in\{1,2\}$, or else $g$ behaves like one of $\id, \mix\id\rev, \mix\rev\id, \mix\rev\rev,\sw, \sw\circ\mix\rev\rev,$ $\sw\circ\mix\id\rev,$ or $\sw\circ\mix\rev\id$ on $X$. We may thus assume the latter holds, for both $X$ and $Y$. Note that each of $h^4, k^4, (h\circ k)^4,$ and $(k\circ h)^4$ behaves like $\id$ on $D$. 

Let $S,T\subseteq D$ be such that $S\oo T$. Then there is a self-embedding of $\Pi$ which sends $S$ and $T$ into $X$ and $Y$, respectively. To see this, let $H$ be a $\oo$-downward closed subset of $D$ without a $\oo$-largest element which contains $S$ and is disjoint from $T$. Let $(\Pi, H)$ denote the structure obtained from $\Pi$ by adding the unary relation $H$ to the language. Let $\Pi\restriction_{X\cup Y}$ be the structure induced by $X\cup Y$, and let $(\Pi\restriction_{X\cup Y}, X)$ denote the structure obtained from $\Pi\restriction_{X\cup Y}$ by adding the unary relation $X$ to the language. Then the structures $(\Pi, H)$ and $(\Pi\restriction_{X\cup Y}, X)$ are isomorphic, and any isomorphism from the first to the latter is an embedding of $\Pi$ as desired.

If $g[X]\oo g[Y]$, then by the above there exists an embedding $\iota$ of $\Pi$ such that $\iota\circ g[X]\subseteq X$ and $\iota\circ g[Y]\subseteq Y$. Let $g'=(\iota\circ g)^4$. Then since $h^4$ and $k^4$ behave like $\id$ on $D$, $g'$ is a function mon-generated by $g$ which behaves like $\id$ on $X$ and $Y$, and such that $\up(g'[X],g'[Y])$ holds.

If on the other hand $g[Y]\oo g[X]$, then pick an embedding $\iota$ of $\Pi$ such that $\iota\circ g[Y]\subseteq X$ and $\iota\circ g[X]\subseteq Y$. Let $g'=(i\circ g)^8$. Since $(h\circ k)^4$ and $(k\circ h)^4$ both behave like $\id$ on $D$, $g'$ is a function mon-generated by $g$ which behaves like $\id$ on $X$ and $Y$, and such that $\up(g'[X],g'[Y])$ holds, thus proving our claim.

Take any $a_1,\ldots,a_n\in D$ with $a_i\oo a_{i+1}$ for all $1\leq i\leq n-1$. By homogeneity of $\Pi$, for each $1\leq i\leq n$ there exists $\alpha\in\Aut(\Pi)$ such that $\alpha(a_j)\in X$ for all $1\leq j\leq i$ and $\alpha(a_j)\in Y$ for all $i+1\leq j\leq n$. Then $g'\circ\alpha$ diagonalizes the sets $\{a_1,\ldots,a_i\}$ and $\{a_{i+1},\ldots,a_n\}$ while preserving $\oo$ on $\{a_1,\ldots,a_n\}$. By repeated applications of such functions, we obtain a function $f$ mon-generated by $g$ such that $\up(f(a_i),f(a_{i+1}))$ holds for all $1\leq i\leq n-1$. It then follows by homogeneity of $\Pi$ and topological closure that $g$ mon-generates a canonical function from $\Pi$ to $\Pi$ whose image is a diagonal set, and Proposition~\ref{prop:canonical} implies that $\G$ contains $\Aut(D;<_i)$ for some $i\in\{1,2\}$.
\end{proof}

\begin{lem}\label{lem:id orbits}
Let $\G$ be a closed supergroup of $\Aut(\Pi)$. Let $c_1,\ldots, c_n\in D$ be distinct and let $g$ be a function mon-generated by $\G$ which is canonical as a function from $(\Pi, c_1,\ldots, c_n)$ to $\Pi$. Suppose $g$ behaves like $\id$ on some infinite orbit of $\Aut(\Pi,c_1,\ldots,c_n)$. Then $g$ behaves like $\id$ on all infinite orbits of $\Aut(\Pi,c_1,\ldots, c_n)$, or else $\Aut(D;<_i)\subseteq \G$, for some $i\in\{1,2\}$.
\end{lem}
\begin{proof} 
Let $X$ be an infinite orbit of $\Aut(\Pi,c_1,\ldots,c_n)$ on which $g$ behaves like $\id$. Let $Y$ be an orbit in the same column as $X$; without loss of generality, suppose $X\ot Y$. Then by Lemmas~\ref{lem:orbits} and~\ref{lem:diagonalize}, either $\Aut(D;<_i)\subseteq \G$, for some $i\in\{1,2\}$, or $g$ preserves $\oo$ on $X\cup Y$. We may thus assume the latter holds. Therefore, $g$ either behaves like $\id$ or $\mix\id\rev$ on $Y$, and like $\id$ or $\mix\id\rev$ between $X$ and $Y$. 

We claim that if $g$ behaves like $\mix\id\rev$ on $Y$, then $\G$ contains $\Aut(D;<_i)$ for some $i\in\{1,2\}$. To see this, observe first that in that situation, $g$ mon-generates $\mix\id\rev$. Therefore, if $g$ behaves like $\mix\id\rev$ between $X$ and $Y$, then we can replace $g$ by $\mix\id\rev\circ g$, which behaves like $\id$ on $Y$, like $\mix\id\rev$ on $X$, and like $\id$ between $X$ and $Y$. Hence, in any case we may assume that $g$ behaves like $\id$ between $X$ and $Y$, like $\id$ on one of $X$ or $Y$, and like $\mix\id\rev$ on the other. Without loss of generality, suppose $g$ behaves like $\id$ on $X$ and $\mix\id\rev$ on $Y$. Now take any $a_1,\ldots,a_m\in D$ with $a_i\ot a_{i+1}$ for all $1\leq i\leq m-1$. By homogeneity of $\Pi$, for any $1\leq i\leq m$ there exists $\alpha\in\Aut(\Pi)$ such that $\alpha(a_j)\in X$ for all $1\leq j\leq i$ and $\alpha(a_j)\in Y$ for all $i+1\leq j\leq m$. Then $g\circ\alpha$ preserves $\ot$ on $\{a_1,\ldots,a_i\}$ and reverses $\ot$ on $\{a_{i+1},\ldots,a_m\}$, while preserving $\oo$ on $\{a_1,\ldots,a_m\}$. By repeated applications of such functions, we can change the order of $\{a_1,\ldots,a_m\}$ arbitrarily with respect to $\ot$ while preserving $\oo$. Therefore $g$ mon-generates a function $f$ such that $\up(f(a_i),f(a_{i+1}))$ for all $1\leq i\leq m-1$. It then follows by homogeneity of $\Pi$ and topological closure that $g$ mon-generates a canonical function $h\colon \Pi\To\Pi$ whose image is a diagonal set, and Proposition~\ref{prop:canonical} implies that $\G$ contains $\Aut(D;<_i)$ for some $i\in\{1,2\}$.

Now let $Z$ be an infinite orbit in the same row as $X$. By the same argument as above, either $g$ behaves like $\id$ on $Z$ or $\Aut(D;<_i)\subseteq \G$ for some $i\in\{1,2\}$. It follows that $g$ behaves like $\id$ on all infinite orbits of $\Aut(\Pi,c_1,\ldots, c_n)$, or else $\Aut(D;<_i)\subseteq \G$, for some $i\in\{1,2\}$.
\end{proof}

\begin{lem}\label{lem:rigid}
Let $\G$ be a closed supergroup of $\Aut(\Pi)$. Let $c_1,\ldots,c_n\in D$ be distinct and $g\colon (\Pi,c_1,\ldots,c_n)\To \Pi$ be a canonical function mon-generated by $\G$. Suppose that $g$ behaves like $\id$ on and between all infinite orbits of $\Aut(\Pi,c_1,\ldots,c_n)$. Then $g$ behaves like $\id$ everywhere, or else $\Aut(D;<_i)\subseteq \G$, for some $i\in\{1,2\}$.
\end{lem}
\begin{proof}
Assume that $g$ does not behave like $\id$ everywhere, say without loss of generality that $g$ does not preserve $\ot$ on $D$. Since $g$ is canonical as a function from $(\Pi,c_1,\ldots,c_n)$ to $\Pi$ and behaves like $\id$ on and between all infinite orbits of $\Aut(\Pi,c_1,\ldots,c_n)$, there exist $c\in\{c_1,\ldots,c_n\}$ and $0\leq j\leq n-1$ such that the rows $R_j, R_{j+1}$ satisfy $R_j\ot\{c\}\ot R_{j+1}$ and either $\{g(c)\}\ot g[R_j]$ or $g[R_{j+1}]\ot\{g(c)\}$. Suppose without loss of generality the latter holds.

Let $0\leq k\leq n-1$ be such that the columns $C_k, C_{k+1}$ satisfy $C_k\oo\{c\}\oo C_{k+1}$, and let $U=(R_j\cup R_{j+1})\cap(C_k\cup C_{k+1})$. Then $g[U]\ot \{g(c)\}$.

First suppose that $g$ diagonalizes $U$ and $\{c\}$; say without loss of generality that $g[U]\oo \{g(c)\}$. Take any $a_1,\ldots, a_m\in D$ with $a_i\oo a_{i+1}$ for all $1\leq i\leq m-1$. By homogeneity of $\Pi$, for each $1\leq i\leq m$ there exists $\alpha\in\Aut(\Pi)$ such that $\alpha(a_i)=c$ and $\alpha(a_l)\in U$ for all $l\neq i$. Then $g\circ\alpha$ diagonalizes the sets $\{a_1,\ldots,a_{i-1},a_{i+1},\ldots,a_m\}$ and $\{a_i\}$, while behaving like $\id$ on $\{a_1,\ldots,a_{i-1},a_{i+1},\ldots,a_m\}$. By repeated applications of such functions, we obtain a function $f$ mon-generated by $g$ such that $\up(f(a_i),f(a_{i+1}))$ holds for all $1\leq i\leq m-1$. It follows from the homogeneity of $\Pi$ and topological closure that $g$ mon-generates a canonical function from $\Pi$ to $\Pi$ whose image is a diagonal set. Therefore, by Proposition~\ref{prop:canonical}, $\Aut(D;<_i)\subseteq\G$ for some $i\in\{1,2\}$.

It remains to consider the case where $g[C_k]\oo \{g(c)\}\oo g[C_{k+1}]$. Then $\Aut(D;\oo)\subseteq\G$, since by a similar argument as above, we can change the order of the elements of any finite subset of $D$ with respect to $<_2$ whilst keeping their order with respect to $<_1$ by repeated applications of functions in $\{g\}\cup\Aut(\Pi)$.
\end{proof}

\begin{lem}\label{lem:permute} Let $\G$ be a closed supergroup of $\Aut(\Pi)$. Let $c_1,\ldots,c_n\in D$ be distinct and $g\colon (\Pi,c_1,\ldots,c_n)\To \Pi$ be a canonical function mon-generated by $\G$. If $g$ behaves like $\id$ on some infinite orbit of $\Aut(\Pi,c_1,\ldots,c_n)$, then either $g$ is mon-generated by $\la\mix\id{t},\mix{t}\id\ra$ or $\Aut(D;<_i)\subseteq \G$ for some $i\in\{1,2\}$.
\end{lem}
\begin{proof} By Lemma~\ref{lem:id orbits}, either $\Aut(D;<_i)\subseteq \G$ for some $i\in\{1,2\}$ or $g$ behaves like $\id$ on all infinite orbits of $\Aut(\Pi,c_1,\ldots,c_n)$. We may thus assume the latter holds. Then by Lemma~\ref{lem:orbits}, either $g$ diagonalizes two infinite orbits in nondiagonal position, in which case $\Aut(D;<_i)\subseteq \G$ for some $i\in\{1,2\}$ by Lemma~\ref{lem:diagonalize}, or $g$ behaves like one of $\id$ or $\mix\id\rev$ between infinite orbits in the same column, and like $\id$ or $\mix\rev\id$ between infinite orbits in the same row. Again, we may assume the latter holds. 

Suppose that there exist infinite orbits $X\ot Y\ot Z$ in the same column of $(\Pi,c_1,\ldots,c_n)$ such that $\neg\cyc_2(g[X],g[Y],g[Z])$ holds. Suppose $g[X]\ot g[Z]\ot g[Y]$. The other cases are proved similarly. 

Take any $a_1,\ldots,a_m\in D$ with $a_i\ot a_{i+1}$ for all $1\leq i\leq m-1$. By homogeneity of $\Pi$, for each $1\leq j\leq m$ there exists $\alpha\in\Aut(\Pi)$ such that $\alpha(a_i)\in X$ for all $1\leq i\leq j-1$, $\alpha(a_j)\in Y$, and $\alpha(a_i)\in Z$ for all $j+1\leq i\leq m$. Then $g\circ\alpha$ reverses $\ot$ between $\{a_j\}$ and $\{a_{j+1},\ldots,a_m\}$, while preserving $\ot$ on $\{a_1,\ldots,a_{j-1},a_{j+1},\ldots,a_m\}$ and preserving $\oo$ on $\{a_1,\ldots,a_m\}$. By repeated applications of such functions, we can change the order of $\{a_1,\ldots,a_m\}$ arbitrarily with respect to $\ot$ while preserving $\oo$. Thus, $g$ mon-generates a function $f$ such that $\up(f(a_i),f(a_{i+1}))$ holds for all $1\leq i\leq m-1$. It follows from homogeneity of $\Pi$ and topological closure that $g$ mon-generates a canonical function from $\Pi$ to $\Pi$ whose image is a diagonal set. Hence, by Proposition~\ref{prop:canonical}, $\Aut(D;<_i)\subseteq\G$ for some $i\in\{1,2\}$.

By the same argument, if there exist infinite orbits $X\oo Y\oo Z$ in the same row of $(\Pi,c_1,\ldots,c_n)$ such that $\neg\cyc_1(g[X],g[Y],g[Z])$ holds, then $\Aut(D;<_i)\subseteq\G$ for some $i\in\{1,2\}$. 

We may henceforth assume that for $i=1,2$ and for all infinite orbits $X<_iY<_iZ$, $\cyc_i(g[X],g[Y],g[Z])$ holds. Suppose there exists $0\leq j\leq n-1$ such that $g$ reverses $\oo$ between the columns $C_j$ and $C_{j+1}$. Then since $\cyc_1(g[X],g[Y],g[Z])$ holds for all infinite orbits with $X\oo Y \oo Z$, it follows that $g[C_{j+1}]\oo\cdots\oo g[C_n]\oo g[C_0]\oo\cdots\oo g[C_j]$. Furthermore, $g$ mon-generates $\mix{t}\id$: for any finite tuple $a$ of elements of $D$, there exists an
embedding $\iota$ of $\Pi$ into the structure induced by $Y\cup Z$ such that $g\circ \iota$ sends $a$ to a tuple of equal type in $\Pi$ as does $\mix{t}\id$; and then we can refer to homogeneity of $\Pi$ to see that $g$ mon-generates a function which agrees with $\mix{t}\id$ on $a$.

By definition, there are $H,H'\subseteq D$ such that $H\oo H'$, $H\cup H'=D$, and such that $\mix{t}\id[H']\oo \mix{t}\id [H]$. We may assume that $g[C_n]\subseteq H$ and $g[C_0]\subseteq H'$. Then $\mix{t}\id\circ g$ is a function mon-generated by $g$ which is canonical from $(\Pi,c_1,\ldots,c_n)$ to $\Pi$ and which preserves $\oo$ on and between all infinite orbits of $\Aut(\Pi,c_1,\ldots,c_n)$. By a similar argument, we can undo a possible shuffling of rows by applying $\mix\id{t}$ if necessary, obtaining a canonical function $h\colon (\Pi,c_1,\ldots,c_n)\To\Pi$ which behaves like $\id$ on and between all infinite orbits of $\Aut(\Pi,c_1,\ldots,c_n)$. By Lemma~\ref{lem:rigid}, either $h$ behaves like $\id$ everywhere or $\Aut(D;<_i)\subseteq\G$ for some $i\in\{1,2\}$. We may thus assume the former holds. That means that by composing $g$ with functions in $\{\mix\id{t},\mix{t}\id\}$ from the left, we have obtained a self-embedding $h$ of $\Pi$. But then $g$ is itself a composite of $h$ with functions in $\{\mix\id{t},\mix{t}\id\}$, and so it is mon-generated by $\la\mix\id{t},\mix{t}\id\ra$.
\end{proof}

\subsection{Proof of Theorem~\ref{thm:JI}.}\label{subsection:proof}

We are now ready to prove Theorem~\ref{thm:JI}, from which it follows that every closed group properly containing $\Aut(\Pi)$ is the join of elements of $\JI$.

\begin{proof}[Proof of Theorem~\ref{thm:JI}]
Let $c_1,\ldots, c_n\in D$ be so that no function in $\G$ agrees with $f$ on $\{c_1,\ldots,c_n\}$. Then by Proposition~\ref{prop:magical}, there is an injection $g\colon D\To D$ mon-generated by $f$ which is canonical as a function from $(\Pi,c_1,\ldots, c_n)$ to $\Pi$ and which agrees with $f$ on $\{c_1,\ldots,c_n\}$. 

Case 1: Suppose $\Aut(D;<_i)\subseteq \G$ for some $i\in\{1,2\}$. Without loss of generality, suppose $i=1$. By Theorem~\ref{thm:cameron5}, $\G$ and $\la\{f\}\cup \G\ra$ are among the following groups:

\begin{itemize}
\item $\Aut(D;\oo)$;
\item $\la\Aut(D;\oo)\cup\{\mix\rev\id\}\ra$;
\item $\la\Aut(D;\oo)\cup\{\mix{t}\id\}\ra$;
\item $\la\Aut(D;\oo)\cup\{\mix\rev\id,\mix {t}\id\}\ra$.
\end{itemize} The result follows, since $\G\subsetneq\la\{f\}\cup \G\ra$.

Case 2: Suppose $\Aut(D;<_i)\nsubseteq \G$ for $i=1,2$ and let $X$ be an infinite orbit of $\Aut(\Pi,c_1,\ldots,c_n)$. Then by Proposition~\ref{prop:canonical}, $g$ behaves like some function $h\in\{\id, \mix\id\rev, \mix\rev\id, \mix\rev\rev, \sw, \sw\circ\mix\rev\rev,$ $\sw\circ\mix\id\rev,\sw\circ\mix\rev\id\}$ on $X$. If $h\notin \G$, then we are done. Otherwise, $h^{-1}\circ g$ is a function which is not mon-generated by $\G$, is canonical as a function from $(\Pi,c_1,\ldots,c_n)$ to $\Pi$, and behaves like $\id$ on $X$. So, by replacing $g$ with $h^{-1}\circ g$, we may assume that $g$ behaves like $\id$ on $X$. 

Therefore, by Lemma~\ref{lem:permute}, $g$ is mon-generated by $\la\mix\id{t},\mix{t}\id\ra$. Since no function in $\G$ agrees with $g$ on $\{c_1,\ldots,c_n\}$, $g$ is not mon-generated by $\G$. Therefore, $\la\mix\id{t},\mix{t}\id\ra\nsubseteq\G$. 

By Lemma~\ref{lem:rigid}, $g$ does not behave like $\id$ between all infinite orbits of $\Aut(\Pi,c_1,\ldots,c_n)$. Therefore, by Lemma~\ref{lem:orbits}, $g$ either behaves like $\mix\id\rev$ between two infinite orbits in the same column or like $\mix\rev\id$ between two infinite orbits in the same row. Suppose the latter holds. Then $g$ mon-generates $\mix{t}\id$. If $\mix{t}\id\notin\G$, then we are done. Otherwise, $\la\mix{t}\id\ra\subseteq\G$. In this case, since $g$ is not mon-generated by $\G$ but is mon-generated by $\la\mix\id{t},\mix{t}\id\ra$, $g$ must also behave like $\mix\id\rev$ between two infinite orbits in the same column. Therefore, $g$ mon-generates $\mix\id{t}$. Since $\la\mix\id{t},\mix{t}\id\ra\nsubseteq\G$, it follows that $\mix\id{t}\in\la\{f\}\cup \G\ra\setminus\G$.
\end{proof}

\section{The 39 closed supergroups of $\Aut(\Pi)$}\label{section:lattice}

We will now determine the precise number of closed supergroups of $\Aut(\Pi)$, proving Theorem~\ref{thm:L},  and compare our result with the estimate from~\cite{CameronPermutations}.

\subsection{Proof of Theorem~\ref{thm:L}}
In order to see that there are at most 39 closed supergroups of $\Aut(\Pi)$, we need the following easy to verify results about the behaviors of compositions of the functions which generate the groups in $\JI$.

\begin{lem}\label{lem:sw} Let $\G$ be a closed supergroup of $\Aut(\Pi)$. Suppose $\sw\in \G$. Then
\begin{itemize}
\item If $\G$ contains one of $\mix\id\rev$ or $\mix\rev\id$, then it contains both $\mix\id\rev$ and $\mix\rev\id$.
\item If $\G$ contains one of $\mix\id{t}$ or $\mix{t}\id$, then it contains both $\mix\id{t}$ and $\mix{t}\id$.
\item If $\sw\circ\mix\rev\rev\in \G$, then $\mix\rev\rev\in \G$.
\item If $\sw\circ\mix\id\rev\in \G$, then $\mix\id\rev\in \G$.
\item If $\Aut(D;<_i)\subseteq \G$ for some $i\in\{1,2\}$, then $\G=\Sym(D)$.
\end{itemize}
\end{lem}

\begin{lem}\label{lem:swrev} Let $\G$ be a closed supergroup of $\Aut(\Pi)$. Suppose $\sw\circ\mix\rev\rev\in \G$. Then
\begin{itemize}
\item If $\G$ contains $\mix\id\rev,\mix\rev\id,\mix\rev\rev,$ or $\sw\circ\mix\id\rev$, then it contains $\sw$.
\item If $\G$ contains one of $\mix\id{t}$ or $\mix{t}\id$, then it contains both $\mix\id{t}$ and $\mix{t}\id$.
\item If $\Aut(D;<_i)\subseteq \G$ for some $i\in\{1,2\}$, then $\G=\Sym(D)$.
\end{itemize}
\end{lem}

\begin{lem}\label{lem:revrev} Let $\G$ be a closed supergroup of $\Aut(\Pi)$. Suppose $\mix\rev\rev\in \G$. Then if $\G$ contains one of $\mix\id\rev$ or $\mix\rev\id$, then it contains both $\mix\id\rev$ and $\mix\rev\id$.
\end{lem}

\begin{lem}\label{lem:swidrev} Let $\G$ be a closed supergroup of $\Aut(\Pi)$. If $\sw\circ\mix\id\rev\in \G$, then 
\begin{itemize}
\item $\mix\rev\rev\in \G$.
\item If $\G$ contains one of $\mix\id{t}$ or $\mix{t}\id$, then it contains both $\mix\id{t}$ and $\mix{t}\id$.
\item If $\Aut(D;<_i)\subseteq \G$ for some $i\in\{1,2\}$, then $\G=\Sym(D)$.
\end{itemize}
\end{lem}

\begin{prop} There are at most 39 closed supergroups of $\Aut(\Pi)$.
\end{prop}
\begin{proof} 
There are at most 25 groups which arise as joins of groups in $\{\la\mix\id\rev\ra,\la\mix\id{t}\ra,\la\mix\rev\id\ra,$ $\la\mix{t}\id\ra,\Aut(D;\oo), \Aut(D;\ot)\}$: by Theorem~\ref{thm:cameron5}, there are 5 closed supergroups of $\Aut(\Pi)$ which contain $\Aut(D;\oo)$ and 4 additional groups containing $\Aut(D;\ot)$; the remaining 16 groups are all possible joins of groups in $\{\la\mix\id\rev\ra,\la\mix\id{t}\ra,\la\mix\rev\id\ra,\la\mix{t}\id\ra\}$. We remark that these are precisely the intersections of closed supergroups of $\Aut(D;\oo)$ with closed supergroups of $\Aut(D;\ot)$, as follows from Figure~\ref{table:pres}.

The remaining closed groups must contain one of $\sw,\sw\circ\mix\rev\rev,$ or $\mix\rev\rev$. By Lemma~\ref{lem:sw}, there are at most 6 additional closed groups containing $\la\sw\ra$:
\begin{enumerate}
\setcounter{enumi}{25}
\item $\la\sw\ra$;
\item $\la\sw,\mix\rev\rev\ra$;
\item $\la\sw,\mix\id\rev\ra$;
\item $\la\sw,\mix\id{t}\ra$;
\item $\la\sw,\mix\rev\rev,\mix\id{t}\ra$;
\item $\la\sw,\mix\id\rev,\mix\id{t}\ra$.
\end{enumerate} By Lemma~\ref{lem:swrev}, there are at most 2 additional groups containing $\la\sw\circ\mix\rev\rev\ra$:
\begin{enumerate}
\setcounter{enumi}{31}
\item $\la\sw\circ\mix\rev\rev\ra$;
\item $\la\sw\circ\mix\rev\rev,\mix\id{t}\ra$;
\end{enumerate} Last, by lemmas~\ref{lem:revrev} and~\ref{lem:swidrev}, there are at most 6 more groups containing $\la\mix\rev\rev\ra$:
\begin{enumerate}
\setcounter{enumi}{33}
\item $\la\mix\rev\rev\ra$;
\item $\la\sw\circ\mix\id\rev\ra$;
\item $\la\mix\rev\rev,\mix\id{t}\ra$;
\item $\la\mix\rev\rev,\mix{t}\id\ra$;
\item $\la\sw\circ\mix\id\rev,\mix\id{t}\ra$;
\item $\la\mix\rev\rev,\mix\id{t},\mix{t}\id\ra$.
\end{enumerate}
\end{proof}

To see that these groups all are distinct, we need to define some more relations on $D$. 
\begin{itemize}
\item $\textsf{R}_1(x,y,z)\Leftrightarrow (\up(x,y)\wedge\dow(y,z)\wedge\up(x,z))\vee (\dow(x,y)\wedge\up(z,y)\wedge\dow(x,z))\\
	\indent\indent\indent\indent\indent\vee (\up(y,x)\wedge\dow(z,y)\wedge\up(z,x))\vee (\dow(y,x)\wedge\up(y,z)\wedge\dow(z,x));$
\item $\textsf{R}_2(x,y,z)\Leftrightarrow \btw_1(x,y,z)\vee\btw_2(x,y,z);$
\item $\textsf{R}_3(x,y,z)\Leftrightarrow \cyc_1(x,y,z)\vee\cyc_2(x,y,z);$
\item $\textsf{R}_4(x,y,z)\Leftrightarrow \sepa_1(x,y,z)\vee\sepa_2(x,y,z);$
\item $\textsf{R}_5(x,y,z) \Leftrightarrow ((x\ot y\ot z)\wedge \cyc_1(x,y,z))\vee ((z\ot y\ot x)\wedge \cyc_1(z,y,x));$
\item $\textsf{R}_6(x,y,x) \Leftrightarrow ((x\oo y\oo z)\wedge \cyc_2(x,y,z))\vee ((z\oo y\oo x)\wedge \cyc_2(z,y,x));$
\item $\textsf{R}_7(x,y,z)\Leftrightarrow ((\cyc_1(x,y,z)\wedge\cyc_2(x,y,z))\vee ((\cyc_1(z,y,x)\wedge\cyc_2(z,y,x));$
\item $\textsf{R}_8(x,y,z) \Leftrightarrow \cyc_1(x,y,z)\wedge\neg\cyc_2(x,y,z);$
\item $\textsf{R}_9(x,y,w,z)\Leftrightarrow (\cyc_1(x,y,w)\wedge\cyc_2(x,y,w)\wedge\cyc_1(y,w,z)\wedge\neg\cyc_2(y,w,z))\\
	\indent\indent\indent\indent\indent\indent\vee(\neg\cyc_1(x,y,w)\wedge\cyc_2(x,y,w)\wedge\cyc_1(y,w,z)\wedge\cyc_2(y,w,z))\\
	\indent\indent\indent\indent\indent\indent\vee(\neg\cyc_1(x,y,w)\wedge\neg\cyc_2(x,y,w)\wedge\neg\cyc_1(y,w,z)\wedge\cyc_2(y,w,z))\\
	\indent\indent\indent\indent\indent\indent\vee(\cyc_1(x,y,w)\wedge\neg\cyc_2(x,y,w)\wedge\neg\cyc_1(y,w,z)\wedge\neg\cyc_2(y,w,z)).$
\end{itemize}

Figure~\ref{table:pres} shows which relations are preserved by the 39 groups listed above. Checking the table is left to the reader. Since no two groups preserve the same subset of relations, Theorem~\ref{thm:L} follows.

\subsection{Discussion}\label{subsec:discussion}
In~\cite{CameronPermutations}, Cameron listed 37 closed supergroups of $\Aut(\Pi)$. His count included the 25 groups which arise as intersections of closed supergroups of $\Aut(D;\oo)$ and $\Aut(D;\ot)$. He then observed that $\la\sw,\mix\id\rev\ra$ behaves like a dihedral group of order 8, with 10 subgroups, 4 of which were already counted in the first 25. This gives 6 additional groups contained in $\la\sw,\mix\id\rev\ra$. By the same argument, he counted 6 additional groups contained in $\la\sw,\mix\id{t}\ra$. 

We discovered that while $\la\sw,\mix\id\rev\ra$ behaves like a dihedral group of order 8, $\la\sw,\mix\id{t}\ra$ does not. It only has 4 proper, nontrivial subgroups: $\la\sw\ra, \la\mix\id{t}\ra, \la\mix{t}\id\ra,$ and $\la\mix\id{t},\mix{t}\id\ra$. Thus there is some asymmetry in the roles of the permutations $\mix\id\rev$ and $\mix\id{t}$. There is a closed group consisting of all permutations which either preserve or reverse both orders of $\Pi$ simultaneously, namely $\la\mix\rev\rev\ra$, but there is no simultaneous action of turns since $\la\mix\id{t}\circ\mix{t}\id\ra=\la\mix\id{t},\mix{t}\id\ra$. Hence, 4 groups counted in~\cite{CameronPermutations} actually coincided with others listed there. On the other hand, some joins of elements in $\JI$ were missing in~\cite{CameronPermutations}: $\la\mix\rev\rev,\mix\id{t}\ra, \la\mix\rev\rev,\mix{t}\id\ra, \la\mix\rev\rev,\mix\id{t},\mix{t}\id\ra,$ $\la\sw,\mix\rev\rev,\mix\id{t}\ra,$ $\la\sw\circ\mix\id\rev,\mix\id{t}\ra$, and $\la\sw,\mix\id\rev,\mix\id{t}\ra$.

\begin{figure}[h]\tiny
\begin{tabular}{|c|c|c|c|c|c|c|c|c|c|c|c|c|c|c|c|c|c|c|c|c|}\hline
   & $\oo$ & $\btw_1$ & $\cyc_1$ & $\sepa_1$ & $\ot$ & $\btw_2$ & $\cyc_2$ & $\sepa_2$ & $\st$ & $\up$ & $\dow$ & $\textsf{R}_1$ & $ \textsf{R}_2$ & $\textsf{R}_3 $ & $\textsf{R}_4 $ & $\textsf{R}_5 $ & $\textsf{R}_6 $ & $\textsf{R}_7 $ & $\textsf{R}_8$& $\textsf{R}_9$\\\hline
a    &  x  &  x  &  x  &  x  &     &  x  &     &  x  &     &     &     &     &	x  & 	 &  x  &     &     &     &     &   \\\hline
b    &  x  &  x  &  x  &  x  &     &     &  x  &  x  &     &     &     &     &	   &  x  &  x  &     &  x  &  x  &  x  &  x\\\hline
c    &     &  x  &     &  x  &  x  &  x  &  x  &  x  &     &     &     &     &	x  & 	 &  x  &     &     &     &     &   \\\hline
d    &     &     &  x  &  x  &  x  &  x  &  x  &  x  &     &     &     &     &     &  x  &  x  &  x  &     &  x  &  x  &  x\\\hline
e    &     &  x  &     &  x  &     &  x  &     &  x  &  x  &     &     &  x  &  x  &     &  x  &  x  &  x  &  x  &     &  x\\\hline
f    &     &     &     &     &     &     &     &     &  x  &  x  &     &     &	x  &  x  &  x  &     &     &  x  &     &   \\\hline
g    &     &     &     &     &     &     &     &     &  x  &     &  x  &     &	x  &     &  x  &     &     &  x  &  x  &   \\\hline
h    &     &     &     &     &     &     &     &     &     &     &     &  x  &	x  &     &  x  &     &     &     &     &  x\\\hline
i    &  x  &  x  &  x  &  x  &     &     &     &     &     &     &     &     &     &     &     &     &     &     &     &   \\\hline 
j    &     &     &     &     &  x  &  x  &  x  &  x  &     &     &     &     &     &     &     &     &     &     &     &   \\\hline

ab   &  x  &  x  &  x  &  x  &     &     &     &  x  &     &     &     &     &	   &     &  x  &     &     &     &     &   \\\hline
ac   &     &  x  &     &  x  &     &  x  &     &  x  &     &     &     &     &	x  &	 &  x  &     &     &     &     &   \\\hline 
ad   &     &     &  x  &  x  &     &  x  &     &  x  &     &     &     &     &     &	 &  x  &     &     &     &     &   \\\hline 
af   &     &     &     &     &     &     &     &     &     &     &     &     &	x  &	 &  x  &     &     &     &     &   \\\hline 
aj   &     &     &     &     &     &  x  &     &  x  &     &     &     &     &	   &	 &     &     &     &     &     &   \\\hline
bc   &     &  x  &     &  x  &     &     &  x  &  x  &     &     &     &     &	   &	 &  x  &     &     &     &     &   \\\hline 
bd   &     &     &  x  &  x  &     &     &  x  &  x  &     &     &     &     &	   &  x  &  x  &     &     &  x  &  x  &  x\\\hline 
be   &     &  x  &     &  x  &     &     &     &  x  &     &     &     &     &	   &     &  x  &     &  x  &  x  &     &  x\\\hline 
bf   &     &     &     &     &     &     &     &     &     &     &     &     &	   &  x  &  x  &     &     &  x  &     &   \\\hline 
bg   &     &     &     &     &     &     &     &     &     &     &     &     &	   &	 &  x  &     &     &  x  &  x  &   \\\hline 
bh   &     &     &     &     &     &     &     &     &     &     &     &     &	   &	 &  x  &     &     &     &     &  x\\\hline 
bj   &     &     &     &     &     &     &  x  &  x  &     &     &     &     &	   &	 &     &     &     &     &     &   \\\hline 
cd   &     &     &     &  x  &  x  &  x  &  x  &  x  &     &     &     &     &	   &	 &  x  &     &     &     &     &   \\\hline 
ci   &     &  x  &     &  x  &     &     &     &     &     &     &     &     &	   &	 &     &     &     &     &     &   \\\hline 
de   &     &     &     &  x  &     &  x  &     &  x  &     &     &     &  x  &     &	 &  x  &  x  &     &  x  &     &  x\\\hline 
di   &     &     &  x  &  x  &     &     &     &     &     &     &     &     &	   &	 &     &     &     &     &     &   \\\hline 
ef   &     &     &     &     &     &     &     &     &  x  &     &     &     &	x  &	 &  x  &     &     &  x  &     &   \\\hline

abc  &     &  x  &     &  x  &     &     &     &  x  &     &     &     &     &	   &	 &  x  &     &     &     &     &   \\\hline
abd  &     &     &  x  &  x  &     &     &     &  x  &     &     &     &     &	   &	 &  x  &     &     &     &     &   \\\hline  
abf  &     &     &     &     &     &     &     &     &     &     &     &     &	   &	 &  x  &     &     &     &     &   \\\hline  
abj  &     &     &     &     &     &     &     &  x  &     &     &     &     &	   &	 &     &     &     &     &     &   \\\hline 
acd  &     &     &     &  x  &     &  x  &     &  x  &     &     &     &     &	   &	 &  x  &     &     &     &     &   \\\hline
bcd  &     &     &     &  x  &     &     &  x  &  x  &     &     &     &     &	   &	 &  x  &     &     &     &     &   \\\hline
bde  &     &     &     &  x  &     &     &     &  x  &     &     &     &     &	   & 	 &  x  &     &     &  x  &     &  x\\\hline
bef  &     &     &     &     &     &     &     &     &     &     &     &     & 	   &	 &  x  &     &     &  x  &     &   \\\hline
cdi  &     &     &     &  x  &     &     &     &     &     &     &     &     &	   & 	 &     &     &     &     &     &   \\\hline

abcd &     &     &     &  x  &     &     &     &  x  &     &     &     &     & 	   &	 &  x  &     &     &     &     &   \\\hline
\end{tabular}
\caption{Preservation table}\label{table:pres}
\end{figure}

\bibliographystyle{alpha}
\bibliography{permutation.bib}

\newcommand{\etalchar}[1]{$^{#1}$}
\begin{thebibliography}{PPP{\etalchar{+}}11}

\bibitem[Ben97]{Bennett-thesis}
James~H. Bennett.
\newblock {\em The reducts of some infinite homogeneous graphs and
  tournaments}.
\newblock PhD thesis, Rutgers university, 1997.

\bibitem[BF13]{BoettcherFoniok}
Julia B\"ottcher and Jan Foniok.
\newblock Ramsey properties of permutations.
\newblock {\em Electronic Journal of Combinatorics}, 20(1), 2013.

\bibitem[Bod]{Bod-New-Ramsey-classes}
Manuel Bodirsky.
\newblock {New Ramsey} classes from old.
\newblock {\em Electronic Journal of Combinatorics}.
\newblock To appear. Preprint arXiv:1204.3258.

\bibitem[BP]{RandomMinOps}
Manuel Bodirsky and Michael Pinsker.
\newblock Minimal functions on the random graph.
\newblock {\em Israel Journal of Mathematics}.
\newblock To appear. Preprint arXiv.org/abs/1003.4030.

\bibitem[BP11a]{BP-reductsRamsey}
Manuel Bodirsky and Michael Pinsker.
\newblock Reducts of {R}amsey structures.
\newblock {\em AMS Contemporary Mathematics, vol. 558 (Model Theoretic Methods
  in Finite Combinatorics)}, pages 489--519, 2011.

\bibitem[BP11b]{BodPin-Schaefer}
Manuel Bodirsky and Michael Pinsker.
\newblock Schaefer's theorem for graphs.
\newblock In {\em Proceedings of the Annual Symposium on Theory of Computing
  (STOC)}, pages 655--664, 2011.
\newblock Preprint of the long version available at arxiv.org/abs/1011.2894.

\bibitem[BPP13]{42}
Manuel Bodirsky, Michael Pinsker, and Andr\'{a}s Pongr\'acz.
\newblock The 42 reducts of the random ordered graph.
\newblock Preprint arXiv:1309.2165, 2013.

\bibitem[BPT13]{BPT-decidability-of-definability}
Manuel Bodirsky, Michael Pinsker, and Todor Tsankov.
\newblock Decidability of definability.
\newblock {\em Journal of Symbolic Logic}, 78(4):1036--1054, 2013.
\newblock A conference version appeared in the Proceedings of LICS 2011, pages
  321--328.

\bibitem[Cam76]{Cameron5}
Peter~J. Cameron.
\newblock Transitivity of permutation groups on unordered sets.
\newblock {\em Mathematische Zeitschrift}, 148:127--139, 1976.

\bibitem[Cam02]{CameronPermutations}
Peter~J. Cameron.
\newblock Homogeneous permutations.
\newblock {\em Electronic Journal of Combinatorics}, 9(2), 2002.

\bibitem[Hod97]{Hodges}
Wilfrid Hodges.
\newblock {\em A shorter model theory}.
\newblock Cambridge University Press, Cambridge, 1997.

\bibitem[JZ08]{JunkerZiegler}
Markus Junker and Martin Ziegler.
\newblock The 116 reducts of $(\mathbb{Q},<,a)$.
\newblock {\em Journal of Symbolic Logic}, 74(3):861--884, 2008.

\bibitem[Pon11]{Andras-thesis}
Andr\'as Pongr\'acz.
\newblock Reducts of the {H}enson graphs with a constant.
\newblock Preprint, 2011.

\bibitem[PPP{\etalchar{+}}11]{Poset-Reducts}
P\'{e}ter~P\'{a}l Pach, Michael Pinsker, Gabriella Pluh\'{a}r, Andr\'{a}s
  Pongr\'{a}cz, and Csaba Szab\'{o}.
\newblock Reducts of the random partial order.
\newblock Preprint arXiv:1111.7109, 2011.

\bibitem[Sok]{Sok-Directed-graphs-and}
Miodrag Soki\'{c}.
\newblock Directed graphs and {Boron} trees.
\newblock Preprint available from
  http://www.its.caltech.edu/\~{}msokic/SAP3.pdf.

\bibitem[Sok10]{Sokic}
Miodrag Soki\'{c}.
\newblock {\em Ramsey property of posets and related structures}.
\newblock PhD thesis, University of Toronto, 2010.

\bibitem[Tho91]{RandomReducts}
Simon Thomas.
\newblock Reducts of the random graph.
\newblock {\em Journal of Symbolic Logic}, 56(1):176--181, 1991.

\bibitem[Tho96]{Thomas96}
Simon Thomas.
\newblock Reducts of random hypergraphs.
\newblock {\em Annals of Pure and Applied Logic}, 80(2):165--193, 1996.

\end{thebibliography}

\end{document}